\title{Colouring powers of cycles from random lists}
\author {Michael Krivelevich\thanks{
Department of Mathematics, Raymond and Beverly Sackler Faculty of Exact
Sciences, Tel Aviv University, Tel Aviv 69978, Israel. E-mail:
krivelev@post.tau.ac.il. Research supported in part by a USA-Israel
BSF Grant and by a grant from the Israel Science Foundation.}
\and Asaf Nachmias\thanks{
Department of Mathematics, Raymond and Beverly Sackler Faculty of Exact
Sciences, Tel Aviv University, Tel Aviv 69978, Israel. E-mail:
asafnach@post.tau.ac.il.}}

\documentclass[11pt]{article}
\usepackage{latexsym}
\usepackage{amsfonts}
\usepackage{amsmath}
\newtheorem{theo}{Theorem}[section]
\newtheorem{prop}[theo]{Proposition}

\date{}
\begin{document}
\maketitle
\begin{abstract}
Let $C_n^k$ be the $k$-th power of a cycle on $n$ vertices
(i.e. the vertices of $C_n^k$ are those of the $n$-cycle, and two
vertices are connected by an edge if their distance along the cycle is
at most $k$).  For each vertex draw
uniformly at random a subset of size $c$ from a base set $\mathcal{S}$
of size $s=s(n)$. In this paper we solve the problem of determining
the asymptotic probability of the existence of a proper colouring from the lists for all
fixed values of $c,k$, and growing $n$.

\end{abstract}
\section{Introduction}

Let $G$ be a simple and undirected graph. Assign to each vertex
$x$ of $G$ a set $L(x)$ of colours (positive integers). Such an
assignment $L$ of sets to vertices in $G$ is referred to as a
\em{colour scheme} \rm for $G$. An \em L-colouring \rm of $G$ is a
mapping $f$ of $V(G)$ into the set of colours such that $f(x) \in
L(x)$ for all $x\in V(G)$ and $f(x)\neq f(y)$ for each $(x,y) \in
E(G)$. If $G$ admits an $L$-colouring, then $G$ is said to be \em
L-colourable \rm. In case of $L(x)=\{1,\ldots ,k\}$ for all $x \in
V(G)$, we also use the terms \em k-colouring \rm and \em
k-colourable \rm respectively. A graph $G$ is called \em
k-choosable \rm if $G$ is $L$-colourable for every colour scheme
$L$ of $G$ satisfying $|L(x)|=k$ for all $x \in V(G)$. The \em
chromatic number $\chi(G)$ \rm (\em choice number $ch(G)$ \rm) of
$G$ is the least integer $k$ such that $G$ is $k$-colourable
($k$-choosable).

The choosability concept was introduced, independently by Vizing
\cite{Viz} and by Erd\H{o}s, Rubin and Taylor \cite{ERT}.

Denote by $C_n^k$ the $k$-th power of a cycle on $n$ vertices,
i.e. $V(C_n^k)=\{0,\ldots,n-1\}$ and $(i,j) \in E(G)$ if $(i-j)$
mod $n$ $\in \{-k,\ldots,-1,1,\ldots,k\}$. Equivalently, the
vertices of $C_n^k$ are those of the $n$-cycle, and two vertices
are connected by an edge if their distance along the cycle is at
most $k$. Let $\mathcal{S}$ be a set of colours of size $s(n)$.
For each vertex of $C_n^k$ draw uniformly at random $c$ colours
from $\mathcal{S}$ to form a colour scheme $L(c,k,s)$. Denote by
$p(n,c,k,s(n))$ the probability that $C_n^k$ is
$L(c,k,s)$-colourable. The problem of determining what is the
minimum $s(n)$ such that $p(n,c,k,s(n)) = 1-o(1)$ was posed to us
by Maurice Cochand. The problem originated in chemical industry
and it is related to scheduling problems occurring in the
production of colorants. The goal is to get an estimate of the
resources required for the production process. Basically, the
vertices represent jobs, the colours correspond to processors that
perform the jobs, and the edges represent technological
restrictions on the processors' assignments.

Of course, the same problem of colourability from randomly chosen
lists can be posed for other asymptotic families such as complete
graphs $K_n$ and complete bipartite graphs $K_{m,n}$ as well. We
plan to return to this question later.

In this paper we investigate the asymptotic behaviour of $p(n)=p(n,c,k,s)$ for
every fixed $c,k$ and $n$ tending to infinity. We show that if $c \leq k$,
$$ p(n) = \left \{
\begin{array}{ll}
o(1), & s(n)=o(n^{1/c^2}) , \\
e^{-\binom{k}{c}t^{-c^2}(c!)^c} + o(1), & s(n) \sim tn^{1/c^2} , \\
1-o(1), & s(n) = \omega(n^{1/c^2}) , \\
\end{array} \right .
$$

\noindent On the other hand, if $c=k+1$,
$$ p(n) = \left \{
\begin{array}{ll}
o(1), & s(n) \leq c , \\
1-o(1), & s(n) >c , \\
\end{array} \right .
$$

\noindent and if $c > k+1$
$$ p(n) = \left \{
\begin{array}{ll}
o(1), & s(n) < c ,\\
1-o(1), & s(n) \geq c .\\
\end{array} \right .
$$

We shall use the standard asymptotic notation and assumption. In
particular we assume that the parameter $n$ is large enough
whenever necessary. For two functions $f(n)$ and $g(n)$, we write
$f=o(g)$ if $\lim_{n\rightarrow\infty} f/g=0$, and $f=\omega(g)$
if $g=o(f)$. Also, $f=O(g)$ if there exists an absolute constant
$c>0$ such that $f(n)<cg(n)$ for all large enough $n$;
$f=\Theta(g)$ if both $f=O(g)$ and $g=O(f)$ hold; and $f\sim g$ if
$\lim_{n\rightarrow\infty} f/g=1$.

\section{Case $c \leq k$}

For all results in this section we assume $c \leq k$. We shall
see that with probability tending to 1 the graph becomes $L$-colourable
at the same time cliques of size $c+1$ with the same list drawn for
every vertex, vanish.

\begin{prop}
If $s(n) = o(n^{1/c^2})$, then $p(n)=o(1)$.
\end{prop}

\noindent{\bf Proof:} Every set $\{v_1,...,v_{c+1}\}$ of $c+1$
vertices of $C_n^k$ satisfies:
$$ Pr[v_1, ..., v_{c+1} \textrm{ have the same list }]=\binom{s}{c}^{-c};$$
Partition the first $(c+1)\lfloor \frac{n}{c+1}  \rfloor$ vertices of $V(C_n^k)$ \
into $\lfloor \frac{n}{c+1} \rfloor$ disjoint sets of $c+1$ consecutive
(in the sense of the cycle) vertices. Clearly, the events ``The set gets
the same list drawn for all its vertices" are mutually independent.
 Therefore, the probability that
none of those sets has the same list is:
$$ \left [ 1 - \frac{1}{\binom{s}{c}^{c}} \right ]^{\lfloor n/(c+1) \rfloor}
= o(1)\ . $$
Hence with probability tending to $1$, there exists a set of $c+1$ consecutive
vertices, which forms a clique in $C_n^k$, with the same list drawn for every
vertex. Therefore, $C_n^k$ cannot be coloured.

\hfill $\Box$

\begin{theo}\label{th1}
If $s(n) = \omega (n^{1/c^2})$, then $p(n)=1-o(1)$.
\end{theo}
\noindent{\bf Proof:} Fix with foresight a constant $d$
satisfying $d>\frac{c^2(c^2+c-1)}{c-1}$. Order the vertices
$v_1, ..., v_n$ according to the cycle, with an arbitrary starting point.
Call a colour scheme $L$ of $C_n^k$ {\em good} if it satisfies the following
conditions:

\begin{enumerate}
\item There are no $c+1$-cliques with the same list drawn for
every vertex of the clique.
\item There exists a family of sets of $k+1$
consecutive vertices (in what follows: $k+1$-sets) such that:
\begin{enumerate}
\item The lists drawn for each of the sets do not intersect with any of the
$2k$ lists of its neighbours.
\item The number of vertices between two such consecutive $k+1$-sets
is no more than $n^{1/d}$.
\end{enumerate}
\item Every set of at most $n^{1/d}$ consecutive vertices, $U$, has the
following property: every subset $X \subset U$, $|X|\geq c+2$,
satisfies $|X| \leq |\{c : \exists x \in X, c\in L(x)\}|$.
\end{enumerate}
We prove the theorem in two steps. First we show that if $L$ is
good, then $C_n^k$ is $L$-colourable. Secondly, we show that
with probability tending to $1$, $L(c,k,s(n))$ is good.

\vspace{1pc}

Assume $L$ is good, then every set of at most $n^{1/d}$ consecutive vertices,
$V$, can be coloured from $L$. To see this construct a bipartite graph, $H_L$ with
sides $V,C$ where $C=\mathcal{S}$ and $(v,c) \in E(H)$ if $c \in
L(v)$. By Condition $3$, $|X|\leq |N(X)|$ for all $X$ of size $|X|
\geq c+2$. Since the degree of every $v \in V$ is $c$, it follows
that $|X|-|N(X)|\leq 1$ for every $X \subset V$. By the defect
version of Hall's theorem (see \cite{West}) there exists a matching in $H$
which saturates all
but at most one vertex of $V$. This matching ensures us a legal choice of
colours for the saturated vertices. Furthermore, it is clear that if
$x \in V$ is the unsaturated vertex, then $x$ is in some $X
\subset V$ where $|X|-|N(X)| = 1$. This can only happen if
$|X|=c+1$, meaning also that all the lists of $X$'s vertices are
identical. By Condition $1$, $X$ is not a clique. Colour two
non-adjacent vertices of $X$ with the same colour, if $x$ is still
not coloured, colour it with the available colour. This clearly
completes the matching and thus completes an $L$-colouring of $V$.

Take a family of $k+1$-sets as in Condition $2$, find an
$L$-colouring of every $k+1$-set and an $L$-colouring of every run
of vertices between consecutive $k+1$-sets. This is clearly
possible by the above argument and this completes an $L$-colouring
of $C_n^k$.

\vspace{1pc}

We show now that with probability tending to $1$, $L(c,k,s(n))$ is
good.

\vspace{1pc}

{\em Condition 1.} Let $X$ be the random variable counting the
number of $c+1$-cliques with the same list drawn for every vertex.
The number of cliques in $C_n^k$ is $n\binom{k}{c}$ (choose the
leftmost vertex of the clique in $n$ possible manners, then choose
$k$ vertices from its $c$ neighbours to the right in
$\binom{k}{c}$ possible manners), so clearly:
$$E[X] = n \binom{k}{c} \binom{s}{c}^{-c} = o(1),$$
So $Pr[X \geq 1] = o(1)$.

\vspace{1pc}

{\em Condition 2.} For a fixed $k+1$-set, the
probability that its lists do not intersect with its neighbours'
lists is clearly more than: \newline
$$\left [\frac{\binom{s-2kc}{c}}{\binom{s}{c}} \right ]^{k+1} = 1 -
\Theta(\frac{1}{s});$$
Hence, the chance that the $k+1$-set's lists are not entirely disjoint from its
neighbours' lists is at most $\Theta(1/s)$. \newline
Partition the first $(3k+3) \lfloor \frac{n}{3k+3} \rfloor$ vertices of $V(C_n^k)$
into disjoint $k+1$-sets, with
distance of $2k+1$ between them (this is a mere technicality to keep the
events independent). So, by the union bound, the event that there exists a run
of at least $n^{1/d}$ consecutive vertices in which every $k+1$-set
(in the former division) has failed to satisfy the list
disjointness condition occurs with probability at most:
$$n [\Theta(1/s)]^{\Theta(n^{1/d})} = o(1);$$

\vspace{1pc}

{\em Condition 3.} We show that for a given set of
$\lfloor n^{1/d} \rfloor$ consecutive vertices, $U$, the probability that
Condition $3$ does not hold is $o(1/n)$, thus with probability
tending to $1$, every such $U$ satisfies Condition $3$. This obviously
implies that, with probability tending to $1$, every set of at
most $n^{1/d}$ vertices satisfies Condition 3.
For a fixed $U$, construct a bipartite graph $H_L$ as before. The probability that
Condition $3$ does not hold is clearly:

$$ Pr[\exists X \subset V, |X|\geq c+2, |X| > |N(X)|] \leq
\sum _{i=c+2} ^{\lfloor n^{1/d} \rfloor} \binom{\lfloor n^{1/d} \rfloor}{i} \binom{s}{i-1}
\left [ \frac{\binom{i-1}{c}}{\binom{s}{c}}\right]^i ,$$

\vspace{1pc}
The first term in the summation implies the choice of such $X$, the
second implies the choice of $N(X)$ and the third is the
probability that the neighbours of $X$ are in $N(X)$, i.e., that
the lists of the vertices of $X$ were chosen from the $i-1$
colours (at most) of $N(X)$. This is $o(\frac{1}{n})$:

\begin{eqnarray*}
\sum _{i=c+2} ^{\lfloor n^{1/d} \rfloor} \binom{\lfloor n^{1/d} \rfloor}{i} \binom{s}{i-1}
\left [ \frac{\binom{i-1}{c}}{\binom{s}{c}}\right]^i &\leq&
\sum _{i=c+2} ^{n^{1/d}} \left ( \frac{en^{1/d}}{i} \right )^i
\left (\frac{es}{i-1} \right)^{i-1}
\left [ \frac {\left(\frac{e(i-1)}{c} \right )^c}
              {\left(\frac{s}{c} \right )^c} \right ]^{i}\\
&\leq& \sum _{i=c+2} ^{n^{1/d}} \left ( \frac {e^{c+2} n^{1/d}(i-1)^{c-2}}{s^{c-1}} \right)^i
\frac{i-1}{es}\\
 &\leq& \frac{n^{1/d}}{es} \sum _{i=c+2} ^{n^{1/d}}
\left ( \frac{e^{c+2}n^{\frac{c-1}{d}}}{s^{c-1}} \right )^i \leq
\frac{n^{1/d}}{es} \frac {\left ( \frac{e^{c+2}n^{\frac{c-1}{d}}}{s^{c-1}} \right )^{c+2}}
{1 - \left ( \frac{e^{c+2}n^{\frac{c-1}{d}}}{s^{c-1}} \right )}\\
&=& e^{(c+2)^2}(1+o(1)) \frac{n^{\frac{c^2+c-1}{d}}}{s^{c^2}s^{c-1}}
\leq e^{(c+2)^2}(1+o(1)) \frac{n^{\frac{c^2+c-1}{d} -
\frac{c-1}{c^2}}}{s^{c^2}}\\
&=& o(\frac{1}{n})\ .
\end{eqnarray*}

\hfill $\Box$

\begin{theo}
For any constant $t>0$, if $s(n) \sim tn^{1/c^2}$, then $\lim _{n
\to \infty} p(n) = e^{-\binom{k}{c}t^{-c^2}(c!)^c}$.
\end{theo}
\noindent{\bf Proof:} The same calculations as before show that if
$s(n) \sim tn^{1/c^2}$, Condition $2$ and Condition $3$ hold with
probability tending to $1$. Thus, the problem reduces to
calculating the asymptotic probability of the appearance
of a $c+1$-clique with identical lists drawn for each vertex.
Using Brun's Sieve (see, e.g., \cite{AS}, Chapter 8) we show that the
number of such $c+1$-cliques has Poisson distribution.

Let $m=n\binom{k}{c}$ be the number of $c+1$-cliques in $C_n^k$,
as explained previously. Let $B_i$, $i=1,\ldots ,m=n\binom{k}{c}$,
be the event that the $i$-th $c+1$-clique has the same list drawn
for it. For disjoint $c+1$-cliques,
 it is clear that the
matching $B_i$'s are independent. Let $X_i$, $i=1,\ldots ,m$, be $B_i$'s indicator random
variables, and let $X = \sum _{i=1,\ldots ,m} X_i$ denote the random variable counting
the number of $c+1$-cliques having the same list drawn for its
vertices. We wish to estimate the probability that $X=0$. Define: \newline

$$\mu = \binom{k}{c}t^{-c^2}(c!)^c,$$
$$ S^{(r)} = \sum _{i_1,\ldots ,i_r} Pr[B_{i_1}\cap \ldots \cap B_{i_r}],$$
where the sum is over all sets $\{i_1,\ldots ,i_r \}\subset \{1,\ldots ,m\}.$ Now: \newline
$$E[X] = S^{(1)} = \sum_{i=1}^{m} Pr[B_i] = n\binom{k}{c}\binom{s}{c}^{-c}
\to \mu ,$$

\noindent For $r\geq 2$, divide the sum $S^{(r)}$ into two parts:
$$ S^{(r)} = \sum _{i_1,\ldots ,i_r}^{*} Pr[B_{i_1}\cap \ldots \cap B_{i_r}]
+ \sum _{i_1,\ldots ,i_r}^{**} Pr[B_{i_1}\cap \ldots \cap B_{i_r}],$$
where the first summation goes over all $r$-tuples of pairwise disjoint
$c+1$-cliques and second goes over all $r$-tuples of $c+1$-cliques which
are not pairwise disjoint. Note that every $c+1$-clique intersects with
a constant number of other $c+1$-cliques, therefore, the number of
terms in the first summation is $\binom{m}{r} - O(n^{r-1})$ and:
$$\sum _{i_1,\ldots ,i_r}^{*} Pr[B_{i_1}\cap \ldots \cap B_{i_r}]
= \left [ \binom{m}{r} - O(n^{r-1})\right]\left[ \binom{s}{c} ^{-cr} \right]
\to \frac{\mu^r}{r!} ;$$
To deal with the second summation, note that for every $1\leq i\leq r-1$
there are $\Theta (n^i)$ manners to choose $r$ $c+1$-cliques such
that the size of the maximal disjoint family is $i$. The
probability that every $c+1$-clique in this maximal disjoint family
has the same list drawn is clearly $\binom{s}{c}^{-ci}$. If $i<r$,
then there exists at least one vertex not in the maximal disjoint
family whose list is identical to one of the lists of the
cliques in the maximal disjoint family. Therefore:

$$\sum _{i_1,\ldots ,i_r}^{**} Pr[B_{i_1}\cap \ldots \cap B_{i_r}]
\leq \sum _{i=1}^{r-1} \Theta (n^i)\frac{1}{\binom{s}{c}^{ci}}
\frac{1}{\binom{s}{c}} = o(1), $$
Thus $S^{(r)} \to \mu^r/r!$, and by Brun's Sieve, $Pr[X=0]\to e^{-\mu}$.
It follows by the above discussion that $p(n)$ is also asymptotic to
$e^{-\mu}$.


\hfill $\Box$

\section{ Case $c > k$}

For all results in this section we assume $c > k$.

A recent result of Prowse and Woodall \cite{PW} states that for all values
of $n$ and $k$, $ch(C_n^k)=\chi(C_n^k)$. It is easy to see that when
$n \geq k(k+1)$, if $k+1$ divides $n$ then $\chi(C_n^k)=k+1$, otherwise
$\chi(C_n^k)=k+2$. Their result makes the following proposition
true:
\begin{prop}\label{woodall}
If $s=s(n)$ satisfies $s \geq c \geq k+2$, then $p(n)=1-o(1)$.
\end{prop}

Furthermore, it is trivial that if $s=c=k+1$, no limit exists for the probability of a legal
colour assignment.
We present an argument settling the case $c=k+1$ and
$s>c$, which also implies Proposition \ref{woodall} without the use of the
result of Prowse and Woodall.

\begin{theo}
\vspace{1pc}\
\begin{enumerate}
\item If $c=k+1$ and $s>c$, then $p(n)=1-o(1)$.
\item If $c>k+1$ and $s\geq c$, then $p(n)=1-o(1)$.
\end{enumerate}
\end{theo}
\noindent{\bf Proof:} If $s=\omega (1)$, then as shown before,
for an arbitrary set of k+1 consecutive vertices, the probability
that its lists intersect with its $2k$ neighbours lists is
$\Theta(1/s)$. So, with probability tending to $1$, a fixed set of k+1
consecutive vertices will have lists disjoint from its neighbours'
lists, which allows us to first colour the k+1-set and then
complete the colouring greedily: colour by the order of the cycle,
for each vertex choose a colour not chosen by its preceding
neighbours.

Consider now the case $s=O(1)$. We prove for $c=k+1$, $s=k+2$.
This implies Theorem 3, since increasing $c$ increases the probability of
a colouring, and in this case, the proof works with higher values of $s$ as
long as it stays bounded. Assume then $\mathcal{S}=\{1,...,k+2\}$.

Given a set of $2[(k+1)^2 + (k+1)]$ consecutive vertices, divide the
set into the first and last k+1-sets, and into $k+1$ pairs of disjoint k+1-sets.
Such a set will be called {\em good} if it satisfies the
following conditions:
\begin{enumerate}
\item Each vertex of the first and the last k+1-set has the list $\{1,...,k+1\}$.
\item For every $1\leq i \leq k+1$, each vertex of the second member of the $i$th pair has
the list $\{1,...,k+1\}$ and each vertex of the first member has the list $\{1,...,k+2\}-\{i\}$.
\end{enumerate}
The probability that a set of consecutive $2[(k+1)^2 + (k+1)]$ vertices
satisfies these conditions is:
$$ \left (\frac{1}{k+2} \right)^{2[(k+1)^2 + (k+1)]}\,.$$
After choosing  $\lfloor \frac{n}{2[(k+1)^2 +(k+1)]}\rfloor$ disjoint
consecutive sets in $V(C_n^k)$ in the usual manner,
the probability that none of them has the above properties tends to $0$ and
therefore a good set exists with probability tending to $1$.

Assume we have a good set, colour $C_n^k$ in the following
manner:\newline
Colour the last k+1-set with $(1,...,k+1)$ by that order.
Continue colouring greedily the vertices after the last k+1-set,
by the order of the cycle (again, for each vertex choose a colour
not chosen by its preceding neighbours) until we complete colouring
the first k+1-set (this can be done since every vertex after the first
k+1-set has exactly $k$ coloured neighbours).
For each $1\leq i \leq k+1$ colour as follows:\newline
To the first $i-1$ vertices of the first member of the $i$th pair,
give colours $(1,...,i-1)$ by that order. Give the $i$th vertex of
that member the colour $k+2$, and colour the rest of the first member
greedily. To the first $i$ vertices of the second member of the
$i$th pair, give colours $(1,...,i)$ by that order. Give the rest
of the vertices the same colours as in the matching vertices of
the first member of the $i$th pair.
It is easy to check that this colouring is possible, and
hence the theorem is proven.

\hfill $\Box$

\noindent{\bf Acknowledgement.} The authors wish to thank Maurice
Cochand for posing the problem and for explaining its
motivation/industrial interpretation.

\end{document}